\documentclass[10pt]{IEEEtran}
\usepackage{amsmath,amssymb, amsfonts,amsthm, multirow}
\usepackage{algorithm,algorithmic,bbm,bm, enumitem, xcolor}
\usepackage[colorlinks=true,citecolor=cyan,linkcolor=blue]{hyperref}
\newtheorem{theorem}{Theorem}

\newtheorem{corollary}[theorem]{Corollary}

\newtheorem{lemma}[theorem]{Lemma}

\newtheorem{DE}{Definition}[section]
\newtheorem{PRO}[DE]{Procedure}
\newtheorem{TEMP}[DE]{Template}
\newcommand{\rerror}{{\rm err}}

\newcommand{\rR}{{\rm R}}
\newcommand{\rT}{{\rm T}}
\newcommand{\rH}{{\rm H}}
\newcommand{\rS}{{\rm S}}

\newcommand{\bE}{{\mathbf{E}}}

\def\cA{{\cal A}}

\def\cE{{\cal E}}

\def\cG{{\cal G}}
\def\cK{{\cal K}}
\def\cL{{\cal L}}

\def\cN{{\cal N}}
\def\cP{{\cal P}}
\def\cR{{\cal R}}

\def\cT{{\cal T}}
\def\int{\mbox{int}}

\usepackage{ifpdf}

\usepackage{cite}

  \usepackage[pdftex]{graphicx}
\usepackage{array}

\usepackage{stfloats}
 \usepackage{url}

\begin{document}
%
\title{Stochastic defense against complex grid attacks}

\author{
\IEEEauthorblockN{Daniel Bienstock and Mauro Escobar}
\IEEEauthorblockA{
Columbia University, NY,
United States\\
\{dano, me2533\}@columbia.edu}
}

\maketitle

\begin{abstract}
We describe stochastic defense mechanisms designed to detect sophisticated grid attacks involving both physical actions (including load modification) and sensor output alteration.  The initial attacks are undetectable under a full AC power flow model even assuming ubiquitous sensor placement, while hiding large line overloads.  The defensive techniques apply network control actions that change voltages in a random fashion, and additionally introduce (random) low-rank corrections to covariance matrices. 
\end{abstract}


\noindent 

\section{Introduction}
Recent events and research efforts have highlighted the potential for powerful coordinated attacks on power grids that combine disruption or modification of sensor data with physical actions.  Such attacks may succeed in hiding from operators undesirable system conditions, long enough that physical damage or automatic shutdown of equipment takes place, an undesirable and potentially risky outcome.

We propose  defensive techniques to be deployed when a high-fidelity attack on a power grid is suspected. The attack is assumed to be \textit{partial} in the sense that only a subset of buses and lines are attacked, but this subset is unknown by the system controller. These techniques involve two ideas:
\begin{itemize}
  \item [(a)] using network resources to randomly change power flow quantities, especially voltages and, in particular 
  \item [(b)] changing the covariance structure of e.g. voltages in a manner       unpredictable by the attacker. The specific version of this idea that we analyze introduces a low-rank adjustment to the covariance of phase angles.
\end{itemize}
A precise definition of the attack is given in Section \ref{sec:idattack}.
 Our defensive techniques focus on the phase following the initial attack, and aim to expose inconsistencies in the modified sensor data stream which is output by the attacker.  

We justify such defenses by pointing out that the possibility of dangerous ``cyber-physical'' attacks of high-fidelity and with sparse signatures has already been indicated in the literature (see discussion below). The data component of the attacks is designed to pass a stringent test, namely that the falsified data satisfies the full AC power flow equations (\cite{andersson1,bergenvittal,overbyebook}) at every bus and line.  The data attack is coordinated with a physical attack encompassing various types (in particular, line tripping, or load modifications as considered in this paper) that results in a dangerous system condition, e.g. a line overload. The data modification hides this overload, with the result that sensor data received by operators is both unimpeachable and portrays safe system operation.
We term these attacks ``ideal'' because, while sparse, they do assume technical sophistication and the ability to coordinate physical action and computation. Sparsity is a goal for the attacker because it increases the likelihood of undetectability long enough for the overload to lead to line tripping (typically several minutes).  Putting aside the actual feasibility of such attacks, the computational challenge is significant.  Prior work (e.g. \cite{sankar16}) has already addressed this problem. As further justification for the study of defenses, we develop a new optimization formulation that quickly computes cyber-physical attacks on large transmission systems with thousands of buses.  We stress, however, that we focus on the defensive mechanisms.

For presentational ease we formally describe the attack structure first, in Sections \ref{sec:idattack} and \ref{sec:follow-up}, with the defensive mechanisms given in Section \ref{sec:defense}.

\subsection{Prior work}
The possibility of cyber- or cyber-physical attacks on power grids has yielded  mathematical
work designed to detect and reconstruct such attacks. See \cite{liureiteretal,poor,columbia2,expose,columbia1,react,sankar16,liang16,deng17,danmolzahn,javadattack,vukovicsoudansandberg,liyilmazwang,kimtong,kimtongthomas,dansandberg,deep,valenzuelaintrusion1,davis,anwar,moslemi,mousavian}.

The starting point of this work is the currently used ``State Estimation'' procedure whereby sensor readings are used together with a linearized model of power flow in order to estimate other system parameters.  In its simplest form, this procedure uses the linearized, or DC,  power flow model \cite{andersson1,bergenvittal,overbyebook} (also see Section \ref{sec:notation})
\begin{align} \label{DCflow}
  & B \theta \ = \ P^g - P^d,
\end{align}  
where $B$ is the bus susceptance matrix, $\theta$ is the vector of phase angles, and $P^g$ and $P^d$ are (respectively) the vectors of active power generation and load.  Sensors, which may not be ubiquitous, report phase angles and statistical estimation procedures can be used to recover missing readings as well as other operational data.  If the estimated data fails to satisfy \eqref{DCflow}, an anomalous condition is construed. As discussed in the above works, an attacker that is able to modify sensor output may be able to alter the true phase angles $\theta$ through a perturbation $\delta$ in the null space of $B$; the vector $\theta' \doteq \theta + \delta$ is thus consistent with the equations \eqref{DCflow}.  The resulting attack is thus considered \textit{undetectable} as per the state estimation criterion. Attacks may include a simultaneous physical component which modifies the underlying network, thus providing an additional challenge to a system operator.
As shown in \cite{liureiteretal} such attacks may be sparse (i.e. $\delta$ has small support); computation of an optimally sparse attack is considered in \cite{poor}.

More sophisticated attacks address the AC model of power flows, which is a more accurate description of the actual system behavior than the DC-based equations \eqref{DCflow} (see Section \ref{sec:notation}); attacks that are undetectable under this scenario stand a better chance of actually succeeding in practice while creating a risky condition.  

When the attacker does not have unhindered access to sensors, or if e.g. the result of the attack is that sensors stop reporting, sophisticated techniques may still be brought to bear in order to identify, for example, the topology modification.  See \cite{columbia1, columbia2, expose, react}. Under the assumption that the physical attack disconnects lines, that the defender knows the data-attacked zone of the grid, and structural assumptions about the attack, the techniques in \cite{expose} recover the attack.

With regards to attack detection, some the above work relies on  DC-based state estimation; the model in \cite{expose} uses the AC power flows model. \cite{deng17} studies combined physical- and data-attacks together with countermeasures to detect intrusions, such as reliance on trusted sensors and tracking of power system equivalent impedance; tests in the IEEE 9-bus, 14-bus, 30-bus, 118-bus, and 300-bus test power systems are shown.  Attacks that modify admittances are considered in \cite{danmolzahn}, which also uses inconsistencies in AC the current-voltage relationship to pinpoint the location of an attack. \cite{deng17} relies on estimated impedance computations to localize an attack. \cite{valenzuelaintrusion1} performs PCA (principal component analysis) on the covariance of power flows to discover anomalies by inspecting changes in the smaller eigenvalue modes, also see \cite{lexiepmu,anwar,moslemi,mousavian,pmuuk,wangzhangzhang}.

Of importance to this paper, \cite{sankar16} (also see \cite{liang16}) has shown that at least
in principle, AC-undetectable cyber-physical attacks that result on overloads are
possible. \cite{sankar16} describes a bilevel (two-phase) optimization procedure that
computes a line to be disconnected by the attacker and corresponding data modifications so that a resulting line overload is hidden under the AC power flow model; attacks on the IEEE 24-bus model are reported.  These results add
justification to the study of intelligent defense strategies.

\section{Notation}\label{sec:notation}
We represent AC power flows using the polar representation. The voltage at a
bus $k$ is of the form $V_k = |V_k| e^{j \theta_k}$ where $j = \sqrt{-1}$. A line
$km$ is described by using the standard ``$\pi$'' model which includes series impedance, line charging and transformer attributes. See, e.g. \cite{overbyebook,bergenvittal}. Under this model,
the complex currents injected into line $km$ at bus $k$ and $m$ (respectively)
are given by the formula
\begin{align}\label{picurrent}
  & \left(^{I_{km}}_{I_{mk}}\right) \ = \ Y_{km} \left(^{V_k}_{V_m}\right),
\end{align}
where $Y_{km}$ is the branch admittance matrix for line $km$; 
the complex power injected into line $km$ at $k$ equals
$p_{km} + j q_{km} \ = \ V_k I^*_{km}.$
Here $p_{km} = p_{km}(V_k, V_m)$ and $q_{km} = q_{km}(V_k, V_m)$ are real-valued quadratic functions of the voltages
at $k$ and $m$, which can be summarized in the form
\begin{align}\label{abbreviate}
  & p_{km} + j q_{km} \ = \ S_{km}(|V_k|, |V_m|, \theta_k, \theta_m).
\end{align}
The complex power
flow and angle limits on a line $km$ are denoted by $S_{km}^{max}$ and
$\theta_{km}^{max}$ (respectively), the voltage limits at a bus $k$ are
given by $V_k^{min}$ and $V_k^{max}$, and the active and reactive limits
at a generator bus $k$ are indicated by
$P_k^{g,min}, \ P_k^{g,max}$ and  $Q_k^{g,min},  \ Q_k^{g,max}$ (respectively). 

Given a bus $k$ we denote by $\delta(k)$ the set of all lines of the form $km$.  $\cN$ is the set of buses (we write $n = |\cN|$) and $\cG$ is the set of generator buses\footnote{For simplicity we assume at most one generator per bus.}; given a set of buses $S$ we denote by $\partial S$ (the boundary of $S$) the
subset of buses of $S$ that are incident with a line with an end not in $S$.

We model AGC (Automatic Generation Control) as follows.
There is a selected subset of generators $\cR\subseteq\cG$ (the participating
generators) and 
parameters $\alpha_k \ge 0$ for $k \in \cR$ (the participation factors) with
$\sum_{k \in \cR} \alpha_k = 1$. If aggregate net active power generation changes by some value $\Delta$, with generator $k \in \cR$
changing its output by $\alpha_k \Delta$.

The susceptance matrix $B$ of the DC power flow model \eqref{DCflow} is defined by $B_{kk}=\sum_{km\in\delta(k)}1/x_{km}$ for any bus $k$, $B_{km}=-1/x_{km}$ for any line $km$, and $B_{km}=0$ otherwise; where $x_{km} > 0$ is the reactance of line $km$.  When the underlying network is connected, the solution to system
\eqref{DCflow} has an important attribute, namely that it has one
degree of freedom. Any arbitrary bus 
(the \textit{reference} bus) $r$ can be selected, and 
 $\theta_j$ set to zero.  With this proviso, the system \eqref{DCflow}
has a unique solution,  of the form
\begin{align}\label{refsol}
  & \theta = \breve B_r (P^g - P^d)
\end{align}
where $\breve B_r$ is an appropriate pseudo-inverse of $\hat B$ which depends on the choice of the reference bus.

\vspace{-.2in}
\section{The Initial Phase of an Attack}\label{sec:idattack}
In this section we  discuss, for motivational reasons,
high-fidelity attacks combining physical disruption and data intrusion. We refer the reader to \cite{sankar16} or \cite{danmolzahn} for additional discussions and algorithms.
This paper does not concern the computation of realistic attacks of this nature; however in Section \ref{sec:attackcomp}
we provide a nonlinear, single-phase optimization problem to compute attacks.
We show that the numerical solution to this problem scales well to systems with thousands of buses, with running times in the tens of seconds or less on a standard computer. These experiments
do not prove that the attacks, though sparse, could easily be executed. Rather
they show that the attack computation is tractable, thus 
providing added justification for
studying sophisticated and scalable defense mechanisms.

Attack models considered in prior work allow the attacker different capabilities. Regardless of the model, Template \ref{initatt} given below (similar to one
in \cite{sankar16}) broadly outlines the structure of an AC-undetectable attack.  We use the term ``initial'' to indicate that the attack comprises actions taken at
the start of the attack. Later we will discuss a ``follow-up'' phase that follows the initial attack.

In the template, conditions (a), (f) and (g) amount to a strong form of undetectability. We will provide examples of large scale systems that are susceptible to attacks of the form (a)-(g).  
\begin{center}
  \fbox{
    \begin{minipage}{0.9\linewidth}
\vspace{-.1in}
      \hspace*{.9in} \begin{TEMP}{Initial Attack}\label{initatt} \end{TEMP}
      \begin{itemize}
      \item [(a)] It is assumed that at each bus $k$ there is a sensor measuring voltage at $k$ and current at each line $km \in \delta(k)$.
\item [(b)] The attacker has selected a (sparse) subset $\cA$ of buses, as well as a target line $uv$ within $\cA$ that will be overloaded. 
\item [(c)] The attacker's physical actions are of two types: The attacker may modify \textit{loads} at buses in $\cA$, and may disconnect lines with both ends in $\cA$.
\item [(d)] For any bus $k \in \cA$, the attacker may modify
  data
   provided by a sensor located at $k$. 
\item [(e)] Actions (c)-(d) have to be carefully timed\footnote{See discussion below.}.
\item [(f)] The data received by the control center satisfies complete fidelity as per AC power flow equations and shows all system limits being satisfied, while in actuality line $uv$ is overloaded.
\item [(g)] When the attack includes load changes, secondary response (i.e. AGC response) is taken into
  account by the attacker.
\vspace{-.05in}
\end{itemize}
    \end{minipage}
  }
\end{center}
Note that we allow loads to be modified, but not generation.  In our numerical examples we enforce that $\cG \cap \cA = \emptyset$, out of a
perception that generator sites are more carefully protected.

Point (e) requires some discussion.  Immediately following any physical modification to a system, we can
expect a change in voltages (magnitudes and phase angles) and even to system frequency, the latter especially when net loads are changed. More properly, system dynamics will undergo a change.  Understanding the precise nature of that change is a substantial computational task.  The current state-of-the-art involves a numerical simulation that alternates between simulation of true dynamical behavior at generators (the so-called swing equation) with AC power flow updates.  This combined computation will
typically run much slower than the actual dynamics, and assumes correct knowledge
of the underlying transmission system.  Under adversarial attack that e.g. modifies the topology, the rapid success of such a computational approach to identifying
the current grid state seems uncertain.  And once action (d) (i.e., modification of sensor data) is taken a completely falsified, and consistent view of the system
is being presented.

We next present mathematical conditions that we will impose so
as to guarantee undetectability. \textit{True} data will be the true physical data.
In contrast,
\textit{reported} data is that which is actually received by the control center
and includes the attacker's modifications.
The true data will be given by the (voltage, current) pair of vectors $(V^\rT, I^\rT)$ whereas the reported data will be given by $(V^\rR, I^\rR)$.

An important 
requirement for the reported data is {\bf current-voltage consistency:}
\begin{align}\label{currvolt}
  &\left(^{I^\rR_{km}}_{I^\rR_{mk}}\right) \ = \ Y_{km} \left(^{V^\rR_{k}}_{V^\rR_{m}}\right),
\end{align}
(i.e.) equation \eqref{picurrent}. 
This condition will be enforced in the computation given below in an indirect fashion (also see \cite{danmolzahn} for a different use of this requirement). In general, of course,  an attacker might only seek approximate consistency, using ambient
noise to hide errors. 
Additionally:
\begin{itemize}[leftmargin=.325in]
\item [(s.1)] On a bus $k \notin \cA$ the true and reported data agree (no
  data modification outside $\cA$, by definition).
\item [(s.2)] At a  bus $k \in \partial\cA$ the attacker is constrained by the condition $V^\rR_k = V^\rT_k$.
  This condition is applied to avoid attack detection, given (a) and the second equation in
  \eqref{currvolt} applied to a line $km$ where $m \notin \cA$.
  \item [(s.3)] On buses $k \in \cA - \partial\cA$ we may have $V^\rR_k \neq V^\rT_k$ and on lines with at least one end in  $\cA - \partial\cA$ the true and reported currents may also differ.
\item [(s.4)] The reported voltages and currents must be consistent with meaningful (complex) power injections. Specifically, consider a bus $k$.  Then 
  $\sum_{km \in \delta(k)} V^\rR_k I^{{\rR}*}_{km}$ equals the power injected into the system at bus $k$, according to the reported data.  If $k \notin \cA$ by definition
  (of reported and true data)
  this sum equals $\sum_{km \in \delta(k)} V^\rT_k I^{{\rT}*}_{km}$ which is the true
  power injected by bus $k$.  On the other hand if $k \in \cA$ the sum may differ
  from the true injection at $k$.
  \item [(s.5)] If the attack causes  a net change in the sum of loads, the resulting AGC-mandated
    change in generator output must be taken into account.
\end{itemize}  
We call condition (s.4) {\bf power-injection consistency}, that is:
\begin{align}\label{injcon}
  \text{} & \quad \sum_{km \in \delta(k)} V^\rR_k I^{\rR *}_{km} = \text{net injection at $k$} \quad \forall k,
\end{align}
where ``net injection'' is the reported net injection on buses in $\cA$ and the true net injection for buses in the complement of $\cA$, denoted $\cA^c$.

Subject to these requirements, the attacker seeks to create a
(true) line overload on $uv$ with both ends in $\cA$,
while the reported data shows safe system operation (voltage, angle, and power flow limits are satisfied).  In Section \ref{sec:attackcomp} we will provide
examples of successful attacks on large systems.

\section{The follow-up phase}\label{sec:follow-up}
Above we described actions to be taken by an attacker in order
to initiate an undetectable attack.  In order for the attack to be truly
successful, i.e. by causing an overloaded line to trip, the attack must
remain undetected or at least unreconstructed for a sufficiently long period
of time, possibly on the order of minutes.  This presents a challenge to the
attacker, as the modified sensor data must paint a falsified picture, yet such
sensor data cannot be constant and more generally must follow a realistic
stochastic distribution.  Our defenses will exploit this fact.

Let $t = 0$ denote the time at which the initial attack is completed (start of AGC). In analogy to our notation for the initial problem, at time $t > 0$ we denote by $V_k^\rR(t)$ and $V_k^\rT(t)$ be reported and true voltages at $t$ and similarly with currents. Reported data for  $\cA$ will be manufactured by the attacker 
aiming to approximately satisfy current-voltage consistency \eqref{currvolt} and
power-injection consistency \eqref{injcon}.

In addition, in this work we assume that the attack is perpetrated when
ambient conditions (in particular loads) are, approximately, constant. Let
us denote by $V_k^\rR(0)$ the voltage at a bus $k$ computed by the initial attack, i.e.
$$V_k^\rR(0) \ \doteq \ |V_k^\rR| e^{j \theta^\rR_k}$$ and likewise define the current $I_{km}^R(0)$ on line $km$.  The statement
that ambient conditions are approximately constant, post attack, can be informally rephrased as 
\begin{align}\label{zerodrift}
  & V_k^\rR(t) \, \approx \, V_k^\rR(0)  \ \forall \, k, \text{ and } I_{km}^\rR(t) \, \approx \, I_{km}^\rR(0) \ \forall \, km. \end{align}
If ambient conditions are approximately constant \eqref{zerodrift} will hold (statistically) for any bus $km$ not in the attacked zone $\cA$ but are otherwise a requirement for the attacker.

Two types of attack have been used in the literature.  First, the ``noisy data'' attack in our setting works as follows:
\begin{center}
  \fbox{
    \begin{minipage}{0.9\linewidth}
 \vspace{-.1in}      
      \hspace*{.9in} \begin{TEMP}{Noisy Data Attack}\label{noisydata}\end{TEMP}
      \vspace{-.05in}      
      At time $t > 0$ the attacker reports at each bus $k \in \cA$ a voltage $V_k^\rR(t) \ = \ V_k^\rR(0)  + \bm{\nu_k(t)}.$ \\
      Here $\bm{\nu_k(t)}$ is a random value drawn from a small variance, zero mean distribution\footnote{We use boldface to indicate random variables.}.  Likewise the attacker reports for each line $km$ with both ends in $\cA$, currents
      \begin{align}\label{noisyI1}
        & \left(^{I^\rR_{km}(t)}_{I^\rR_{mk}(t)}\right) \ = \ \left(^{I^\rR_{km}(0)}_{I^\rR_{mk}(0)}\right) + \left(^{\bm{\mu_{km}(t)}}_{\bm{\mu_{mk}(t)}}\right)
        \end{align}
      where $\bm{\mu_{km}(t)}, \bm{\mu_{mk}(t)}$ are drawn from zero mean distributions with small variance.
       \vspace{-.05in}      
    \end{minipage}
  }
\end{center}
Note that these definitions satisfy requirement \eqref{zerodrift}, and approximately satisfy current-voltage consistency.  As a functionally equivalent alternative to \eqref{noisyI1} the attacker could simply set
\begin{align}\label{noisyI2}
  &\left(^{I^\rR_{km}(t)}_{I^\rR_{mk}(t)}\right) \ = \ Y_{km} \left(^{V^\rR_{k}(t)}_{V^\rR_{m}(t)}\right),
\end{align}
our analyses below apply to either form.

A second form of attack that has been considered is the {\bf data replay} attack. Here the attacker supplies a previously observed (or computed) pair of time series $V^\rR(t)$ and $I^\rR(t)$ for buses and lines within the set $\cA$.

\subsubsection{Discussion}\label{sec:attackdiscuss}
The reader may recall that in the initial attack computation we enforced that reported voltages in $\partial \cA$ are exact, i.e. equal to the true voltages. In the follow-up phase this condition is necessarily relaxed by the attacker, though this action carries the risk (to the attacker) that current-voltage consistency will not hold, statistically, for some line $km$ with $k \in \partial\cA$ and $m \notin \cA$. Thus e.g. in the noisy data attack template given above the distributions for the $\bm{\nu_k(t)}$, $\bm{\mu_{km}(t)}$ and $\bm{\mu_{mk}(t)}$ should have sufficiently small \textit{variance} relative to the variance of ambient conditions.  Further requirements on such variances will be discussed in Section \ref{sec:covardef}. In any case, when ambient conditions (e.g. loads) are nearly constant, the noisy-data attack may continue   to approximately satisfy current-voltage and power-injection consistency and  thus remain numerically undetectable. The same holds for the data replay version provided the replayed voltages in $\partial\cA$ closely approximate ambient conditions.

In the next section we present defensive mechanisms that repeatedly change voltages in a way that is unpredictable by the attacker. The key observation is that a substantial change to voltages in $\partial\cA$ will cause the noisy-data attack, applied verbatim as in Template \ref{noisydata}, to fail, because of large current-voltage inconsistencies on lines $km$ with $k \in \partial\cA$ and $m \notin \cA$.  Of course, the template need not be applied verbatim, and in particular the attacker may seek to leverage the possibility of sensor error. We will consider this point in the next section.

In \cite{expose} current-voltage consistency is used in a different setting: (i) the attacked zone $\cA$ is known by the defender, (ii) the attacker only disconnects lines. Under a number of assumptions, in particular that there is a matching between $\cA^C$ and $\cA$ that covers all buses in $\cA$ it is shown that the attack can be accurately recovered.  

\section{Defense}\label{sec:defense}
As discussed above, prior work and our computations in Section \ref{sec:attackex} show that it is possible to compute high-fidelity attacks that disguise dangerous network conditions, even in large, complex transmission systems. Other attacks are also potentially conceivable, e.g. impedance changes, transformer tap changes, etc.

In this section we describe a  generic randomized defense strategy that can deployed when a complex attack is suspected. In the analyses we will assume that the attack impacts a proper subset $\cA$ of the system that is unknown to the control center, as was the case above; though the generic defense strategy applies   under more general attacks as well. We will assume that if there is a topology change the network remains connected.

The strategy can be summarized by the following template:

\begin{center}
  \fbox{
    \begin{minipage}{0.9\linewidth}
    \vspace{-.1in}
    \hspace*{.9in} \begin{PRO}{Random Injection Defense.}\label{randdef} \end{PRO}
    \vspace{-.05in} 
      {\bf Iterate:}\\ 
      \hspace*{.1in}{\bf D1:} Choose, for each $k \in \cG$ a (random) value $\bm{\delta_k}$ such that $\sum_{g \in \cG} \bm{\delta_k} = 0$. Command each generator $k \in \cG$ to change its output to $P^g_k + \bm{\delta_k}$.\\ 
      \hspace*{.1in}{\bf D2:} Following the generation change in step {\bf D1} identify inconsistencies in the observed sensor readings.    
    \end{minipage}
  }
\end{center}
Here, an ``inconsistency'' is an incorrect condition satisfied by the
reported data (such as a violation of a power flow law), or stochastic behavior that
is inconsistent with system-wide behavior understood by the control center.
We will describe several concrete versions of this idea below. See
Procedures \ref{pairsdriven}, \ref{pairsdetect} and \ref{covdetect}.

Each iteration would last several seconds, and statistically significant inconsistencies identified by this scheme are flagged as potential evidence of an attack. Below we will describe several specific implementations of the random ingredient; randomness is used because the attacker cannot anticipate the random injections and thus will not be able to instantaneously update the sensor readings within $\cA$. The above strategy could be AGC-like if only generators $k \in \cR$ (the responding generators) are allowed to have $\bm{\delta_k} \neq 0$ and in general it amounts to
a generator redispatch.  The strategy in Procedure \ref{randdef} is likely
to succeed, in particular against the noisy data or data replay attacks, if the generation changes result in significant voltage changes across the system. Lemma~\ref{change} given below explains why a particular implementation of Procedure \ref{randdef} attains this goal.

We note that there is an existing literature on
using network resources so as to change power flow physics in order to
detect structure or faults. See \cite{probe1,probe2,siddharth1,siddharth2,siddharth3}.
Indeed, even though the description of our random injection defense focuses on power injections, one could also consider other
random probing
strategies that change power flows, such as adjusting transformer settings, controlled line tripping, and the use of other technologies such as storage and solid-state devices.

There are several implementations of the generic strategy. Generally  the
defender wants to make the $|\bm{\delta_k}|$ large because to first order
changes in voltage angles are proportional to $\|\bm{\delta}\|_2$, and a large change in phase
angles is likely to give rise to a
significant current-voltage or power-injection inconsistencies in sensor readings in $\partial\cA$, as discussed above.  This idea
forms the basis for a simple, current-consistency based version of Template \ref{randdef} given in Section \ref{currdef}.

An attacker aware that the random injection defense strategy is applied may try to replace e.g. the noisy data attack
with a more careful manipulation of reported data. For example, the attacker could react to a significant change to voltages in $\partial\cA$ by solving a nonlinear, nonconvex system of inequalities designed
to guarantee approximate current-voltage and power-injection consistency.  In addition, any implied load change within $\cA$ must be very small (or it would contradict observed frequencies).  Finally the attacker would need to perform this
computation very quickly, and repeatedly (because the defense will be applied repeatedly).

This online complex computation could in principle be bypassed by the attacker
by considering changes to readings of voltages at buses in 
$\partial\cA$ only; with the remaining voltages in $\cA$ computed
as in Template \ref{noisydata}.  We will term this the \textit{enhanced} noisy data attack.
We remark that the adversary would still have to maintain AC consistency for lines within
$\cA$, which is nontrivial.  Nevertheless, the ability to adjust readings in
$\partial\cA$ beyond what is prescribed by the noisy data attack may provide some flexibility
for the attacker. However, in Section \ref{currdef} we will show that when the random injection defense causes large-enough voltage changes in $\partial\cA$, the enhanced noisy-data attack fails. See Lemma \ref{hoo}.

A more sophisticated defensive idea,  given in Section \ref{sec:covardef}, changes the \textit{stochastics} of power flow data, in particular voltage covariance, and probes the corresponding properties of the reported data.

Our defensive strategies can be easily adjusted if sensors are not available throughout the system, by restricting the tests we perform to sensorized buses and lines.  Of course, the fewer the sensors the more limited the impact of the defense. Indeed, some interesting work (using the standard, DC-equation state estimation) precisely
seeks to perform system identification post attack when only limited sensor
information is available \cite{columbia1,columbia2,expose,react}.  Note that
the attack problem becomes easier (for the attacker) if sensors are not
widespread.  The attacks computed in Section \eqref{sec:attackex} do assume
sensors at every bus, and yet succeed even in large-scale cases.

 \subsection{Controlling voltages through generation changes}\label{sec:firstorder}
 As discussed above, a  goal of the defense is to produce large voltage angle changes
 in buses in $\cA$, with the intention of revealing inconsistencies in reported data on lines between  $\partial\cA$  and $\cA^c$. The defender, of course, does not know the set $\cA$ and thus it is of interest
 to understand when the voltage at any given bus can be changed by
 appropriately choosing the injections $\delta$.  

 In this section we address the task of changing voltage angles through injections.
 First, we will argue by using the DC power flow approximation \eqref{DCflow}, that
 procedure \ref{randdef} does succeed in changing phase angles (see Lemmas \ref{change} and \ref{assuredchange}).  In Section \ref{sec:numac} we will present experiments under the
 AC power model that verify the DC-based results.  And in Section \ref{currdef} we further argue that the voltage changes are large enough to overcome sensor error.

The defensive strategy that we develop, as a specific implementation of the random injection defense strategy \ref{randdef}, assumes that there is a known set $\cT$ of generator buses that are
  known to be ``trusted'', that is to say, we can assume that data from buses in $\cT$ is
  known to be unmodified. This concept is not new; see \cite{davis,poor,6504815bitar,7921602deng} for related discussions. Without such an assumption the entire suite of signals received by the control centers could be falsified and it is questionable whether any meaningful attack reconstruction can be performed. The following template describes the strategy:    
  \begin{center}
  \fbox{
    \begin{minipage}{0.9\linewidth}
 \vspace{-.1in}
 \hspace*{.9in} \begin{PRO}{Pairs-Driven version of Procedure \ref{randdef}.\label{pairsdriven}} \end{PRO}
 \vspace{-.05in}
 At each execution of step {\bf D1}, select a random pair of generator buses $s$ and $t$, both in $\cT$, as well as random $\bm{\Gamma} > 0$, and use the injections
 \begin{align} \label{twoinj}
  & \bm{\delta_s} = \bm{\Gamma}, \ \bm{\delta_t} = \bm{-\Gamma}, \ \text{and} \ \bm{\delta_k} = 0 \ \forall k \neq s, t.
\end{align}  
\vspace{-.1in}
    \end{minipage}
  }
\end{center}

In an application of this defense, define
  $$ \bm{\hat P^g} \ = \ P^g + \bm{\delta}$$

We denote by $\hat B$ the bus susceptance matrix of the network, post attack. This matrix will be different from the original bus susceptance matrix, $B$,
in case of a topology or susceptance attack; thus the control center does
not know $\hat B$.  Recall that 
as stated above we are assuming that the network remains connected, post attack.

 \begin{lemma} \label{change} Suppose $\hat B \theta = P^g - P^d$,  and
   $\hat B \bm{\hat \theta} = \bm{\hat P^g} - P^d$.  Let $k \neq t$ be a bus such that the post-attack network
   contains a path between $s$ and $k$ that does not include $t$.  Then
\begin{align} \label{phasechange}   
  & \bm{\hat \theta_k} - \bm{\hat \theta_t} \ > \ \theta_k - \theta_t.
\end{align}  
 \end{lemma}
 \noindent \textit{Proof.} Equation \eqref{phasechange} does not change if we subtract from every $\bm{\hat \theta_h}$ any  constant, and likewise with the $\theta_h$.  Thus, without loss of generality $\bm{\hat \theta_t} = \theta_t = 0$. Under this
 assumption \eqref{phasechange} reads:
 \begin{align} \label{phch2}
   & \bm{\hat \theta_k} - \theta_k > 0.
 \end{align}
 Let $M$ be the set of buses $p \neq t$ such that
 \begin{itemize}
 \item [(1)] The network contains a path from $s$ to $p$ that avoids $t$, and
   \item [(2)] Subject to (1), $\bm{\hat \theta_p} - \theta_p$ is \textit{minimum}.
 \end{itemize}
 Aiming for a contradiction, we will assume that
  \begin{align} \label{phch2bad}
    & \bm{\hat \theta_p} - \theta_p \ \le \ 0 \quad \text{for } \ p \in M.
  \end{align}
  Showing that \eqref{phch2bad} is false yields \eqref{phch2}. 
 For any line $km$ define the flow value $f_{km} = (\bm{\hat \theta_k} - \bm{\hat \theta_m} - \theta_k + \theta_m)/x_{km}$.  Since $\hat B (\bm{\hat \theta} - \theta) = \bm{\hat P^g} - P^g$, the flow vector $f$ corresponds (under the DC power flow model) to a power flow with $\bm{\Gamma}$ units of generation at $s$, $\bm{\Gamma}$ units of load at $t$, and zero generation and load elsewhere. Note that for any line $km$, $f_{km} > 0$ iff
 \begin{align}\label{posflow}
 & \bm{\hat \theta_k} - \theta_k > \bm{\hat \theta_m} - \theta_m.
 \end{align}
 This observation implies
 \begin{align}\label{posflow2}
 & \bm{\hat \theta_s} - \theta_s > 0.
 \end{align}
     [To obtain this fact,  decompose the flow vector $f$ into a set of path flows from $s$ to $t$ and telescope \eqref{posflow} along any such path.] 
Pick any $p \in M$ and let $P$ be a path from $s$ to $p$ that avoids $t$.  Say $P = v_0, v_1, \ldots, v_i$ where $v_0 = s$ and $v_i = p$, and let $h$ be smallest such 
 that $v_h \in M$.  By \eqref{posflow2}  $s \notin M$, i.e., $h > 0$.  Then by definition of $h$, 
 $\bm{\hat \theta_{v_{h-1}}} - \theta_{v_{h-1}} > \bm{\hat \theta_{v_{h}}} - \theta_{v_{h}},$
 i.e. $f_{v_{h-1}, v_{h}} > 0$.  But by assumption $v_h \neq t$. So there exists some line $v_{h},m$ such
 that $f_{v_{h}, m} > 0$. Therefore using the assumption $\bm{\hat \theta_k} - \theta_k \le 0$ for all $k \in M$, $v_{h} \in M$, and \eqref{posflow}, 
 \begin{align} \label{baddie} 
   & 0 \ge \bm{\hat \theta_{v_{h}}} - \theta_{v_{h}} >  \bm{\hat \theta_{m}} - \theta_{m}.
\end{align}   
 So $m \neq t$, and as a result  by construction there is a path from $s$ to $m$ that avoids $t$. But
 then \eqref{baddie} contradicts the fact that $v_{h} \in M$.  \qed
 \begin{lemma} \label{assuredchange}  Suppose $k$ is any bus and that there are at least two generators available to implement Procedure \ref{randdef}.  Then a pair $s, t$ satisfying the assumptions of Lemma \ref{change} exists. 
\end{lemma}   
 \noindent \textit{Proof.}  Choose $s \in \cR$ such that $s$ is closest to $k$. \qed

 \noindent Note: if $\bm{\Gamma}$ is chosen negative in \eqref{twoinj}, then instead of
 \eqref{phasechange} we obtain $\bm{\hat \theta_k} - \bm{\hat \theta_t} \, < \, \theta_k - \theta_t$ through essentially the same proof.

 For future reference, we state the following analogue of Lemma \ref{change},
 with similar proof (omitted).
  \begin{lemma} \label{nochange} Let Suppose $\theta$ and
   $\bm{\hat \theta}$ be as in Lemma \ref{change}.  Let $k \neq t$ be a bus such that in the post-attack network
   every path between $s$ and $k$ must include $t$.  Then
\begin{align} \label{phasenochange}   
  & \bm{\hat \theta_k} - \bm{\hat \theta_t} \ = \ \theta_k - \theta_t.
\end{align}  
 \end{lemma}

 To analyze the pairs-driven defense we will express phase angles using
 $t$ as the reference bus; see equation \eqref{refsol}. Thus
 $\hat B \bm{\hat \theta} = \bm{\hat P^g} - P^d$ has a unique solution
 with $\bm{\hat \theta_t} = 0$ of
 the form
\begin{align}\label{affinegamma}
  & \bm{\hat \theta} = \breve B_t (\bm{\hat P^g} - P^d) \ = \ \theta + \breve B_t \bm{\delta}
\end{align}
where $\breve B_t$ is an appropriate pseudo-inverse of $\hat B$.  (In this expression $\theta$ is also written w.r.t. to $t$). As a result of Lemma \ref{change} we have:

\begin{lemma} \label{success1}  Let $k$ be any bus. Then, with high probability
  we will have $\bm{\hat \theta_k} = \theta_k + \beta_k \bm{\Gamma}$ for some value $\beta_k > 0$.
\end{lemma}
\noindent \textit{Proof.} Let $k$ be any bus.  Then with high probability, at multiple iterations of the pairs-driven defense the buses $k, s$ and $t$ will satisfy the conditions of Lemma \ref{change}.  Let as assume that bus $t$ is the reference bus. Then $\bm{\hat \theta_t} = \theta_t = 0$, and thus $\bm{\hat \theta_k} > \theta_k$.  Thus, (by \eqref{affinegamma}) $\bm{\hat \theta_k}$ is as desired. \qed

This results suggests the following detection paradigm:
  \begin{center}
  \fbox{
    \begin{minipage}{0.9\linewidth}
 \vspace{-.1in}
 \hspace*{.9in} \begin{PRO}{Pairs-driven detection criterion}\label{pairsdetect} \end{PRO}
 \vspace{-.05in}
 As Procedure \ref{pairsdriven} iterates, for each bus $k$, estimate the
 correlation between 
 $\bm{\hat \theta_k}$ and $\bm{\Gamma}$. The defensive procedure terminates
 when all these estimates are stable. At that point,  any bus $k$ whose
 correlation coefficient is nonpositive is flagged  as suspicious.
    \end{minipage}
  }
\end{center}

\begin{lemma} \label{success}  With high probability the pairs-driven defense will defeat the noisy data and data replay attacks in the sense that each bus whose data is modified will be flagged, and any bus that is not attacked will not
  be flagged.
\end{lemma}
\noindent \textit{Proof.} Let $k$ be any bus.  Then with high probability by Lemma \ref{success1} the control center expects that $\bm{\hat \theta_k} = \theta_k + \beta_k \bm{\Gamma}$ for some value $\beta_k > 0$ unknown to the control center.  This fact yields the desired result since in either attack case, if $k \in \cA$ the signals produced by the attacker at bus $k$ will not have the stated form.  \qed

Lemmas \ref{change} through \ref{success} assume the DC power flows model, which is only a first-order approximation to the AC model we consider here.  In the next section we perform numerical experiments, under the AC power flows model, of the random injection defense \ref{randdef}.  A separate issue concerns ambient \textit{noise}; we need voltage changes to overcome currently found noise levels in measurements.  This issue is taken up in Section \ref{currdef}.

 \subsubsection{Numerical experiments using AC power flows} \label{sec:numac}
The above discussion concerns DC power flows. In order to investigate how voltages change under injection changes, under AC power flows, we perform a experiments using examples from the Matpower library \cite{matpowerpaper}. For each system we perform ten experiments. In each experiment we compute an AC power flow which is constrained to satisfying the given voltage bounds at all generator buses, but not at load buses, as well as power injection constraints and generator limits, while allowing large injection changes in a random subset of generators.
For a non-generator bus $k$, let $V^b_k$ be its voltage in the base case (i.e. the Matpower case), and let $V_{i,k}$ be its voltage in experiment $i = 1, \ldots, 10$.  Finally, define
$$ \text{score}(k) \doteq \max_{1 \le i \le 10} \frac{|V_{i,k} - V^b_k|}{|V^b_k|}.$$
In Table \ref{tab:voltchange}, ``Min Score'' is the minimum score across all non-generator buses. Thus the table provides experimental verification for substantial AC voltage changes under random generator injections.
\begin{table}[ht]
  \vspace{-.2in}        
  \begin{center}
    \caption{AC Voltage Changes\label{tab:voltchange}}
  \vspace{-.1in}            
    \renewcommand{\arraystretch}{1.2}
    \begin{tabular}{| l | r | r |} 
    Case  & Min Score & Average Score   \\  \hline 
      \texttt{case118} & $11.61\%$ & $32.77\%$   \\ 
      \texttt{case1354pegase} & $7.62\%$ & $51.00\%$  \\
      \texttt{case2746wp} & $5.00\%$ & $10.09\%$ \\ [.5ex] \hline
    \end{tabular}
    \renewcommand{\arraystretch}{1}
  \end{center}
\end{table}

\subsection{Overcoming sensor error, and the current-voltage defense}\label{currdef}
If sensor misestimation (i.e., \textit{error}) is present, a strategy based on Procedure \ref{randdef} may fail to detect data inconsistencies if the random power injections cause
voltage changes that are too small as compared to the error. In order to derive a version of Procedure \ref{randdef} that deals with this issue, we next describe a particular implementation of step {\bf D2} which relies on the current-voltage consistency condition \eqref{currvolt} which takes into account
the possibility of sensor error.  Whereas above a phasor (voltage or current quantity) $\phi$ had a true physical value $\phi^\rT$ and a reported value $\phi^\rR$ (which is the value received by the control center), now we will
have the \textit{sensed} value $\phi^\rS$ which is the value actually produced by the sensor.

Due to sensor error, sensed and true data may differ.   For a phasor $\phi$ define
$ \rerror(\phi) \, \doteq \, \phi^{\rS} - \phi^{\rT}.$
In the PMU setting, the TVE (total vector error) criterion   \cite{pmus3}, \cite{pmus4} guarantees
that
\begin{align}\label{tvedef}
  & | \rerror(\phi) | \, < \, \tau |\phi^{\rT}|,
\end{align}  
where  $0 < \tau < 1$ is a tolerance. Standards enforce $\tau = 1\%$, though experimental testing of PMUs shows far smaller errors \cite{italianspmu}.  From \eqref{tvedef} we obtain
\begin{subequations}
\begin{align}
  & (1 - \tau) | \phi^\rT | \, < \, | \phi^{\rS} | \, < \, (1 + \tau) | \phi^\rT | \label{tve2} \\
  & | \rerror(\phi) | \, < \, \tau(1-\tau)^{-1} | \phi^\rS | \label{tve3}.
\end{align}
\end{subequations}
We will describe two sensor-error-aware voltage-current consistency criteria.  An important point is that the current-voltage consistency condition \eqref{currvolt}, combined with estimations of possible sensor error, yields a nonlinear relationship, and appropriately reformulation of this relationship can render useful benefits. To simplify notation we will drop the ``$(t)$'' from phasors though it should be understood throughout. For a line $km$
write
$$ Y_{km} \ = \ \left(^{Y_{km}^{(1)} \ Y_{km}^{(2)}}_{Y_{km}^{(3)} \ Y_{km}^{(4)}}\right) $$
\noindent {\bf Criterion 1.}  We have that $I_{mk}^\rT = Y_{km}^{(3)} V_k^\rT   + Y_{km}^{(4)} V_m^\rT$.  Write
$Z^{(3)}_{km} \doteq [Y_{km}^{(3)}]^{-1}$.  Hence
\begin{align}
  & V^\rS_k  - Z^{(3)}_{km} (I^\rS_{mk} - Y^{(4)}_{km}V^\rS_{m}) \ = \ \nonumber \\
  & \rerror(V_k) \, - \, Z^{(3)}_{km} (\rerror(I_{mk}) - Y^{(4)}_{km} \rerror(V_m)). \label{basiscrit1}
\end{align}
Using the triangle inequality, we obtain
\begin{align}
  & | V^\rS_k - Z^{(3)}_{km} (I^\rS_{mk} - Y^{(4)}_{km}V^\rS_{m}) | \ <  \nonumber \\ & |\rerror(V_k)| +   |Z^{(3)}_{km}|( \, |\rerror(I^\rS_{mk})| + |Y^{(4)}_{km}||\rerror(V_{m}) | \,) \label{rerrors}
\end{align}  
which yields, using \eqref{tvedef} on the first term and \eqref{tve3} on the
last two terms of \eqref{rerrors},
\begin{align}
  & | V^\rS_k - Z^{(3)}_{km} (I^\rS_{mk} - Y^{(4)}_{km}V^\rS_{m}) | \ \le  \nonumber \\
  & \tau | V^\rT_k| + \frac{\tau |Z^{(3)}_{km}|}{1-\tau}(|I_{mk}^\rS| + |Y^{(4)}_{km}|\,|V_m^\rS |). \nonumber
\end{align}
This last expression equals (using the current-voltage relationship)
\begin{align}
  & \tau |Z^{(3)}_{km} (I^\rT_{mk} - Y^{(4)}_{km}V^\rT_{m}) | + \frac{\tau |Z^{(3)}_{km}|}{1-\tau}(|I_{mk}^\rS| + |Y^{(4)}_{km}|\,|V_m^\rS |), \nonumber
\end{align}
which, using the left inequality in \eqref{tve2} and the triangle inequality, is strictly less than
\begin{align}  
  & \frac{2 \tau |Z^{(3)}_{km}|}{1-\tau}(|I_{mk}^\rS| + |Y^{(4)}_{km}| \, |V_m^\rS |). \nonumber
\end{align}
In summary,
\begin{align}
  & | V^\rS_k - Z^{(3)}_{km} (I^\rS_{mk} - Y^{(4)}_{km}V^\rS_{m}) | \ <  \nonumber \\
  & \frac{2 \tau |Z^{(3)}_{km}|}{1-\tau}(|I_{mk}^\rS| + |Y^{(4)}_{km}| \, |V_m^\rS |). \label{crit1}
\end{align}
Under Criterion 1, if, statistically, the reported phasors $V^\rR_k$, $V^\rR_m$, $I^\rR_{mk}$ fail to satisfy \eqref{crit1} line
$km$ is flagged as suspicious.  A similar analysis concerns $V^\rR_k$, $V^\rR_m$, $I^\rR_{km}$.
\noindent {\bf Remark:} By construction, if $k,m \notin \cA$ then line $km$ will not be flagged.

\noindent {\bf Criterion 2.}  Proceeding as above we have
\begin{align}
  & |I^\rS_{km} - Y^{(1)}_{km} V_k^\rS - Y^{(2)}_{km} V_m^\rS| \, < \,  \nonumber \\
  & \frac{\tau}{1-\tau}( |I_{km}^{\rS}| + |Y^{(1)}_{km}| |V_k^\rS| + |Y^{(2)}_{km}| |V_m^\rS|). \label{crit2}
\end{align}

(and similarly with $I_{mk}$). Note: this criterion can be sharpened  when line $km$ is a pure impedance line (no
transformer).
If the reported phasors do not satisfy \eqref{crit2} then the line is flagged.

\subsubsection{Discussion} Note that a line not attacked will not be flagged, as per the TVE condition. Additional criteria can be developed to handle power-injection consistency.
To analyze the effectiveness of these criteria, we turn to the \textit{enhanced} noisy
data attack discussed in Section \ref{sec:attackdiscuss}. To remind the reader, in
this type of attack the voltage readings in $\partial \cA$ can be arbitrarily adjusted. While this action may create inconsistencies on lines with just one end in $\partial \cA$ the
attacker may be able to ``hide'' such inconsistencies if they are small enough relative
to sensor error.

We next show that Criterion 1 alone can suffice to defeat the enhanced noisy data attack (i.e. uncover inconsistencies) when voltage angles are sufficiently changed under our random injection defense. 

To understand these points, consider a bus $k \in \partial\cA$ such that there is a line $km$ with $m \notin \cA$ and also a line $ka$ where $a \in \cA - \partial\cA$. We study an iteration of the random injection defense which (to simplify notation) we assume begins at time $t = 0$. Consider line $ka$ first. To avoid having line $ak$ flagged, the attacker $t$ will need to manufacture a time series  $V^\rR_k(t)$,  $V^\rR_a(t)$ and $I^\rR_{ak}(t)$
that (statistically) satisfy \eqref{crit1}. But under the noisy data attack, on average $V^\rR_a(t) = V^\rR_a(0)$ and $I^\rR_{ak}(t) = I^\rR_{ak}(0)$.  Hence, to defeat Criterion 1, the attacker needs (on average) that
\begin{align}
  & \frac{2 \tau |Z^{(3)}_{ka}|}{1-\tau}(|I_{ak}^\rR(0)| + |Y^{(4)}_{ka}| \, |V_a^\rR(0) |) \ > \label{crit1app} \\
  & | V^\rR_k(t) - Z^{(3)}_{ka} (I^\rR_{ak}(0) - Y^{(4)}_{ka}V^\rR_{a}(0)) | \, = \, | V^\rR_k(t) - V^\rR_k(0)| \nonumber
\end{align}
Now consider line $km$. Since $m \notin \cA$, $V^\rR_m(t) = V^\rS_m(t)$ and $I^\rR_{mk}(t) = I^\rS_{mk}(t)$.
Also, denote:
\begin{itemize}
  \item $V^\rT_k(*) = $ the true voltage at $k$ at the start of the current iteration of the random
injection defense, i.e. the voltage resulting from the injection changes in step {\bf D1}. Then, 
assuming unbiased sensor errors and zero-mean ambient noise,  $V^\rT_m(*)$ will equal the expectation of $V^\rT_m(t)$ during the iteration.
\item Likewise define the current $I^\rT_{mk}(*)$.
  \end{itemize}
Hence, to defeat Criterion 1, the attacker needs (on average) that
\begin{align}
  & \frac{2 \tau |Z^{(3)}_{km}|}{1-\tau}(|I_{mk}^\rT(*)| + |Y^{(4)}_{km}| \, |V_m^\rT(*)|) \ >
  \label{crit1again} \\
  & | V^\rR_k(t) - Z^{(3)}_{km} (I^\rT_{mk}(*) - Y^{(4)}_{km}V^\rT_{m}(*)) | \, = \, | V^\rR_k(t) - V^\rT_k(*)| \nonumber. 
\end{align}
As a result of these observations we have:
\begin{lemma} \label{hoo} Consider buses $k, a, m$ as described above. Suppose that
  \begin{align}
&    | V^\rT_k(*) - V^\rR_k(0)| > \frac{2 \tau |Z^{(3)}_{ka}|}{1-\tau}(|I_{ak}^\rR(0)| + |Y^{(4)}_{ka}| \, |V_a^\rR(0) |)\nonumber \\ & \hspace{.2in}  + \ \frac{2 \tau|Z^{(3)}_{mk}| }{1 - \tau}(|I_{mk}^\rT(*)| + |Y^{(4)}_{km}| \, |V_m^\rT(*)|). \label{hoo2}
  \end{align}
Then it is impossible for the enhanced noisy data attacker to statistically satisfy Criterion 1 on both lines $ka$ and $km$. 
\end{lemma}
\noindent \textit{Proof.} As argued above, the attacker needs both \eqref{crit1app} and \eqref{crit1again} to hold.  Since
$$| V^\rT_k(*) - V^\rR_k(0)| \le | V^\rR_k(t) - V^\rR_k(0)| + | V^\rR_k(t) - V^\rT_k(*)|,$$ 
summing \eqref{crit1app} and \eqref{crit1again} we obtain a contradiction to \eqref{hoo2}. \qed \\
\noindent {\bf Comment:} This lemma highlights how large changes in voltages caused by the random injection defense challenge the attacker.
  
 \subsubsection{Experiment} Next we describe a set of experiments involving the current-voltage defense applied to the
 attack given in Section \ref{sec:attackex}.  The current defense was implemented as follows:
 \begin{itemize}
 \item For any generator bus $k \notin \cR$, $|\delta_k| \le \epsilon P_k^g$.  We used values $\epsilon = 0.01, 0.05$.
 \item The set of responding generators, $\cR$, was of cardinality $200$. For $k \in \cR$  $|\delta_k|$ can be arbitrarily large. We chose $\delta_k > 0$ with probability $1/2$.
\item No generator may exceed its limits (voltage or generation), but subject to all these conditions we maximize $\sum_{k \in \cG} |\delta_k|$.
\end{itemize}   
 In Table \ref{tab:current-defense}, we perform the above analysis on the lines $(k=1139,a=1137)$
 and $(k=1139,m=1110)$ with $\tau = 0.01$.  ``Ratio'' is the ratio of the left-hand side to the right-hand side of expression \eqref{hoo2}. We see that the
 condition for Lemma \ref{hoo} is amply satisfied. A similar analysis pertains to line $(1141,1361)$, the
 other line connecting $\cA$ to its complement. 
\begin{table}[ht]
  \begin{center}
\vspace{-.1in}    
    \caption{Current-voltage defense.\label{tab:current-defense}}
    \renewcommand{\arraystretch}{1.2}
    \begin{tabular}{| c | r | r |} 
      & Experiment 1 & Experiment 2   \\  \hline 
      $\epsilon$ & 0.01 & 0.05   \\ 
      $\sum_{k\in\cG} \delta_k^+$ & $289.01$ & $964.77$  \\
      $\sum_{k\in\cG} \delta_k^-$ & $174.47$ & $256.04$  \\ [.5ex] \hline 
      \multicolumn{3}{c}{Line $(k=1139,a=1137)$} \\ [.5ex] \hline
      $|V_a^\rR(0)|\angle\theta^\rR_a(0)$    & $1.0919\angle-6.993^\circ$  & $1.0919\angle-6.993^\circ$   \\
      $I_{ak}^\rR(0)$    & $-0.0275 +0.0281j$  & $-0.0275 +0.0281j$ \\ [.5ex] \hline 
      \multicolumn{3}{c}{Line $(k=1139,m=1110)$} \\ \hline
      $|V_m^\rT(*)|\angle\theta^\rT_m(*)$    & $1.0309\angle-7.822^\circ$  &  $1.0391\angle-7.848^\circ$  \\
      $I_{mk}^\rT(*)$    & $0.0905-0.4976j$   & $0.1289-0.4901j$ \\ [.5ex] \hline
\multicolumn{3}{c}{Voltages at $k=1139$} \\ \hline
      $|V_k^\rR(0)|\angle\theta^\rR_k(0)$    & $1.0919\angle-6.991^\circ$  & $1.0919\angle-6.991^\circ$   \\
$|V_k^\rT(*)|\angle\theta^\rT_k(*)$    & $1.0104\angle-7.822^\circ$  & $1.0187\angle-7.936^\circ$  \\ [.5ex] \hline 
\multicolumn{3}{c}{Lemma \ref{hoo} applied to bus $k = 1139$ } \\ \hline
Ratio & $1.913$  & $1.732$\\ \hline
    \end{tabular}
    \renewcommand{\arraystretch}{1}
  \end{center}
\end{table}

\section{Covariance Defense}\label{sec:covardef}
 In this section we describe an elaboration of the pairs-driven defense \ref{pairsdriven}; the elaboration is motivated by the fact that real-world PMU data streams exhibit non-generic stochastic structure in (for example) voltage angles \cite{detective,mynaspi,lexiepmu}. In particular, covariance matrices across several time scales have very low rank (typically smaller than 10).
Our defense will defeat both the noisy-data and data-replay attacks, under
appropriate assumptions.

As before, we assume that the buses
in a certain set $\cT$ are \textit{trusted}.  The emphasis of the methods in
this section is that we aim to modify the \textit{covariance} matrix of phase
angles, whereas the 
random injection defense in Procedures \ref{randdef} or \ref{pairsdriven} change the \textit{average} voltage values.  Such a change
should prove more difficult for the attacker to correctly counteract since
such a correction involves an estimation that requires time, during which
the attacker will be producing incorrect data.

We additionally assume
that the attacker's data stochastics are \textit{stationary} (i.e. the
parameters of the stochastic process do not change as a function of time).
This implies in particular that the attacker does not react to the covariance defense by changing the stochastics. 
Below we will discuss, however, why reacting to the defense would prove very
difficult.    We also assume that ambient conditions are also stationary.

To describe  the defense we need some definitions. For a pair of buses $s, t \in \cT$,  define the vectors $u^{s,t}$ and $v^{s,t}$ by
\begin{subequations}
 \begin{align}
\hspace{-.1in}   & u^{s,t}_{s} = 1,\ u^{s,t}_{t} = -1, \ \text{and} \  u^{s,t}_{k} = 0 \ \forall \, k \neq s, t \label{udef}\\
\hspace{-.1in}   & v^{s,t} \ \doteq \ \breve B_t u^{s,t}. \label{vdef}
 \end{align}
 
\end{subequations}
Formally, the covariance defense works as follows. Let $t_1, t_2$ be
two fixed trusted buses, and  let $\cP$ be a real-valued, zero-mean, positive variance  probability distribution. The defense has two phases.

\noindent {\bf (I)} During an initial phase, post the suspected attack, 
for $i = 1, 2$, we compute the matrix\\
\noindent  $\sigma^2_{\bm{\theta^{\rR,i}}}$ = covariance matrix of observed phase angles, expressed with respect to reference bus $t_i$.\\
\noindent {\bf (II)} After the initial phase, we perform iterations
as in Procedure \ref{pairsdriven}, as follows.   We randomly (uniformly) 
 choose one pair of buses of the form $(s, t)$ where $s \in \cT$ and $t = t_1$
 or $t_2$, and we draw a random value $ \bm{\Gamma}$ from the distribution $\cP$, independently from
 the stochastics of the attacker and ambient stochastics.  We then apply step
 {\bf D1} of Procedure \ref{pairsdriven} using this triple $(s, t, \bm{\Gamma})$.

 Throughout the second phase, we compute the following
 matrix, for $i = 1, 2$,
 \begin{itemize}
 \item $\sigma^2_{\bm{\hat \theta^{\rR,i}}}$ = covariance matrix of observed phase angles with respect to reference bus $t_i$.
 \end{itemize}
 The defense concludes when the estimates for these two matrices become stable.

   \begin{center}
  \fbox{
    \begin{minipage}{0.9\linewidth}
 \vspace{-.1in}
 \hspace*{.9in} \begin{PRO}{Covariance-driven detection criterion}\label{covdetect} \end{PRO}
 \vspace{-.05in}
At termination of the defense, we
 flag a bus $k$ as \textit{suspicious} if, for both $i = 1$ and $2$, the difference between the $(k, k)$ entry of $\sigma^2_{\bm{\hat \theta^{\rR,i}}}$ and  
the corresponding entry of $\sigma^2_{\bm{\theta^{\rR,i}}}$ is smaller than $\lambda$, defined by
\begin{subequations}
\begin{align}
  & \lambda \doteq \frac{\sigma^2_{\bm{\Gamma}}}{|\cT| - 1} \omega, \ \text{where} \label{lambdadef}\\
  & \omega \doteq \min_{s,t,j}\{ (v^{s,t}_j)^2 \, : \, v^{s,t}_j \neq 0 \}. \label{omegadef}
\end{align}
\end{subequations} 
    \end{minipage}
  }
\end{center}

\noindent

This concludes the description of the defense, with analysis given in Lemmas \ref{easy}, \ref{gotit} and  \ref{success}.  

In preparation for those results, suppose that at some point in phase {\bf(II)} the pair $(s, t_i)$ has been selected. Let
 \begin{itemize}
 \item [(a.1)] $\bm{\hat \theta^{\rT,i}}$ is the vector of true voltage phase angles, using $t_i$ as the reference bus.
 \item [(a.2)]  $\bm{\theta^{\rT,i}}$ describes the true vector of voltage phase angles, had the power injections in the defense \textit{not} been applied at that point of time.  It
   is also given using $t_i$ as the reference bus.
 \end{itemize}

Lemma \ref{easy} given next concerns the relationship between these last two vectors.
\begin {lemma} \label{easy} Suppose that at some point in {\bf(II)} the
  pair $(s, t_i)$ is being used. Then $\bm{\hat \theta^{\rT,i}} \, = \, \bm{\theta^{\rT,i}} + \bm{\Gamma} v^{s,t_i}$.
\end{lemma}
\noindent {\em Proof.} Given that the pair $(s, t_i)$ is being used, the
injections in the random defense, per equation \eqref{twoinj} are given
by $\bm{\delta} = \bm{\Gamma} u^{s,t_i}$.  The result follows from equation \eqref{affinegamma}. \qed

The next result presents a key feature of the covariance of phase angles.
Recall that  $\bm{\Gamma}$ is drawn independent of all other stochastics.
\begin{lemma} For $i = 1, 2$,
\begin{align}  
  & \sigma^2_{\bm{\hat \theta^{\rT,i}}} \ = \ \sigma^2_{\bm{\theta^{\rT,i}}} + \frac{\sigma^2_{\bm{\Gamma}}}{|\cT| - 1}\sum_{s \in \cT - t_i} v^{s,t_i} (v^{s,t_i})^\top. \label{covarchange}
\end{align}  
\end{lemma}
\noindent {\em Proof.} We proceed by conditioning on the pair $(s, t_i)$ being selected by the defense. Subject to this conditioning, by Lemma \ref{easy} the covariance of $\bm{\hat \theta^{\rT,i}}$ equals
$$\sigma^2_{\bm{\theta^{\rT,i}}} + \sigma^2_{\bm{\Gamma}} v^i (v^i)^\top \, + \, \text{covar}(\bm{\theta^{\rT,\text{ref}}}, \bm{\Gamma} v^i).$$ The last term in this expression is zero, by the independence assumption on $\bm{\Gamma}$.  The result follows since each pair is chosen with probability $(|\cT| - 1)^{-1}$. \qed

\begin{lemma} \label{gotit}
  Let $k$ be any bus.  Then for at least one of $i = 1$ or $2$, the $(k,k)$ entry
  of $\sigma^2_{\bm{\hat \theta^{\rT,i}}}$ is at least as large  as the
  corresponding entry 
  of $\sigma^2_{\bm{\theta^{\rT,i}}}$, plus $\lambda$ (defined as in \eqref{lambdadef}).
\end{lemma}
\noindent {\em Proof.} Without loss of generality, there is a path between $k$ and $t_1$ that avoids $t_2$.  Note that the pair $(s, t_2)$ with $s = t_1$ is one of the
pairs available for the defense. By Lemma \ref{change}, we have
that $v^{t_1, t_2}_{k} > 0$ and so $v^{t_1, t_2}_{k} \ge \omega^{1/2}$.  Considering equation \eqref{covarchange} for $i = 2$ we see that one of the terms in the sum corresponds to $s = t_1$.  As just
argued, the $(k,k)$ entry of this term is at least $\omega$.  The $(k,k)$ entries in
the remaining terms of the sum are nonnegative (since each term is a positive-semidefinite matrix).  Thus the result follows.\qed

\begin{lemma}\label{finally} The suspicious labels computed by the covariance defense are correct.
\end{lemma}
\noindent {\em Proof.} Consider first a bus $k$ that is not attacked. For such a bus, by
definition, $\bm{\hat \theta^{\rT,i}}_k = \bm{\hat \theta^{\rR,i}}_k$, for  both
$i = 1, 2$.  Thus, by Lemma \ref{gotit}, bus $k$ is not flagged as suspicious. On the other hand, suppose $k$ is attacked. Then under either
the noisy-data or data-replay attacks
the $(k,k)$ entry of $\sigma^2_{\bm{\hat \theta^{\rR,i}}}$ will be equal to the
corresponding entry of $\sigma^2_{\bm{\theta^{\rR,i}}}$, by the stationarity assumption (the attacker does not change stochastics when the defense is implemented).  Hence bus $k$ is flagged. \qed

\noindent {\bf Remarks:} {\bf(1)} Recall \eqref{lambdadef} and \eqref{omegadef}. The
quantity $\omega$ depends on the bus susceptance matrix $B$, only. Hence by
choosing $\sigma^2_{\bm{\Gamma}}$ large enough we can make $\lambda$ large. Also
note that in the above proofs, we can restrict the set $\cT$ to a subset of
size $2$, again helping $\lambda$ attain large values.

\noindent {\bf (2)} Given a pair $(s,t_i)$ used in the defense, by Lemma \ref{change} any entry $v^{s, t_i}_k$ is positive if there is a path from
bus $k$ to $s$ that avoids $t_i$. Let $A$ denote the set of such buses.
By Lemma \ref{nochange}, for $k \notin A$,  $v^{s, t_i}_k = 0$.

Thus, in the term $v^{s,t_i} (v^{s,t_i})^\top$ in \eqref{covarchange}
the entire submatrix with rows and columns in $A$ is positive, and
the remaining entries in $v^{s,t_i} (v^{s,t_i})^\top$ are zero. By adjusting the proof of Lemma \ref{gotit} we conclude that for $k$ and $m$ in $A$, the entry $(k,m)$ of $\sigma^2_{\bm{\hat \theta^{\rT,i}}}$ is at least as large  as the
  corresponding entry 
  of $\sigma^2_{\bm{\theta^{\rT,i}}}$, plus $\lambda$.  Thus a submatrix of the covariance matrix will change via the defense (and not just the diagonal entries).If the network is guaranteed to be $2$-connected \cite{ahuja} one can in fact prove that the entire matrix must change.

The covariance defense has an additional important feature, namely that the sum on
the right-hand side of \eqref{covarchange} has rank-$|\cT|-1$ (as shown next) whereas we expect
the left-hand side of \eqref{covarchange} to have low rank.
\begin{lemma}\label{ranklemma} For $i = 1, 2$ the vectors $v^{s,t_i} = \breve B u^i$ as in \eqref{vdef} are
  linearly independent.  Hence the second term in \eqref{covarchange} has rank
  at least $|\cT| - 1$.
\end{lemma}
\noindent \textit{Proof}.  The $|\cT| - 1$ vectors $u^{s,t_i}$ arising from all pairs $(s, t_i)$
under consideration are linearly independent, by construction in \eqref{udef}.  Hence
the corresponding vectors $v^{s,t_i}$ are also linearly independent. \qed

Lemma \ref{ranklemma} highlights the challenges faced by the attacker, even if the attacker is
aware that the covariance defense is being deployed. The attacker will have to alter
reported data in a way consistent with an appropriate rank change, but the attacker does
not know the pairs $(s_i, t)$ being used (or the distribution $\cP$).  Such ``learning'' would require data observations, i.e. time, during which the attacker is still expected to produce data readings, producing
an error trail.

\section{Appendix: computation of initial attack}\label{sec:attackcomp}
In this section we describe an optimization formulation that computes AC-undetectable attacks that follow the discussion in Section \ref{sec:idattack}; the attacks can involve line disconnection and load modification (which is actually
computed).  Other physical actions (such as impedance changes or transformer tap changes) could be incorporated into the computation. The purpose of this section is to demonstrate that such attacks
are possible on large, realistic system, and that they may be rapidly computed.
Of course, the actual viability of such attacks may depend on the availability of various resources to the attacker, in particular the ability to intrude and remain undetected long enough to learn data on the transmission system.

As input to the computation we have a set $\cA\subset\cN\backslash\cG$
of buses (the target zone), a set of lines $\cL$ to be disconnected, all with both ends in $\cA$, and a line $uv\notin\cL$ with both ends in $\cA$. Write $\cA^C = \cN \setminus \cA$.
Let $(\hat{S}_k^g=\hat{P}_k^g+j\hat{Q}_k^g)_{k\in\cN}$ and $(\hat{S}_k^d=\hat{P}_k^d+j\hat{Q}_k^d)_{k\in\cN}$ be (respectively) the complex power generation and loads at the time of the attack. We assume that the attacker observes
all these quantities. The initial attack problem is given by the following formulation; an explanation of the variables and constraints will be provided below.
\begin{subequations}
  \label{form:attack}
  \begin{flalign}
    & \mbox{Max}\ (p_{uv}^\rT)^2 + (q_{uv}^\rT)^2\\
    & \mbox{s.t.}\nonumber \\
    & \forall k\in\cA^C\cup\partial\cA,\ |V_k^\rT|=|V_k^\rR|, \ \theta_k^\rT=\theta_k^\rR \label{agreeon}\\
    & \forall k\in\cA,\ -(P_k^{d,\rR} + j Q_k^{d,\rR})= \sum_{km\in\delta(k)} (p^\rR_{km} + j q^\rR_{km})\label{rbalwithinA}\\
    & \hspace{.9cm} -(P_k^{d,\rT} + j Q_k^{d,\rT}) = \sum_{km\in\delta(k) \setminus \cL} (p^\rT_{km} + j q^\rT_{km})\label{tbalwithinA}\\
    & \qquad \qquad P_k^{d,\rR}\geq0, \ P_k^{d,\rT}\geq0 \label{nonnegload}\\
    & \forall k\in\cA^C\backslash\cR:\nonumber\\
    & \quad \hat{P}^g_k - \hat{P}^d_k + j (\hat{Q}_k^g - \hat{Q}_k^d) = \sum_{km\in\delta(k)} (p^\rT_{km} + j q_{km}^\rT) \label{baloutsideA}\\
    & \forall k\in\cR:\hspace{1.15cm} P_k^g-\hat{P}_k^g = \alpha_k\Delta \label{secondaryresponse}\\
    & \quad P^g_k - \hat{P}^d_k + j (Q_k^g - \hat{Q}_k^d) = \sum_{km\in\delta(k)} (p^\rT_{km} + j q_{km}^\rT) \\
    & \forall k\in\cG:\nonumber\\
    & \ P_k^{g,min} \leq P_k^g \leq P_k^{g,max}, \ Q_k^{g,min} \leq Q_k^g \leq Q_k^{g,max} \label{genlimits}\\
    & \forall k\in\cN:\  V_k^{min} \leq |V_k^\rT|,|V_k^\rR| \leq V_k^{max}\\
    & \forall \text{ line } km: \nonumber \\
    & \quad | \theta_{k}^\rR - \theta_{m}^\rR| \le \theta_{km}^{max}; \  | \theta_{k}^\rT - \theta_{m}^\rT| \le \theta_{km}^{max} \mbox{ if } km \notin \cL\label{anglelimits} \\    
    & \quad \max\{\, \|(p_{km}^\rR, q_{km}^\rR) \| \, , \, \|(p_{mk}^\rR, q_{mk}^\rR) \| \,\} \leq S_{km}^{max} \label{linelimits}\\    
    & \quad p^\rT_{km} + j q^\rT_{km} = S_{km}(|V_k^\rT|, |V_m^\rT|, \theta_k^\rT, \theta_m^\rT),\ km \notin \cL \label{flow1}\\
    & \quad p^\rT_{mk} + j q^\rT_{mk} = S_{mk}(|V_m^\rT|, |V_k^\rT|, \theta_m^\rT, \theta_k^\rT),\ km \notin \cL\\    
    & \quad p^\rR_{km} + j q^\rR_{km} = S_{km}(|V_k^\rR|, |V_m^\rR|, \theta_k^\rR, \theta_m^\rR)\\
    & \quad p^\rR_{mk} + j q^\rR_{mk} = S_{mk}(|V_m^\rR|, |V_k^\rR|, \theta_m^\rR, \theta_k^\rR) \label{flow4}
\end{flalign}    
\end{subequations}
This formulation uses the following real-valued variables,
where ``T'' indicates true and ``R'', reported:
\begin{itemize}
\item $|V^\rT_k|,\theta^\rT_k, |V^\rR_k|,\theta^\rR_k\ \forall\mbox{ bus }k\in\cN$ (true and reported voltage magnitudes and angles)
\item $P_k^{d,\rT},Q_k^{d,\rT},P_k^{d,\rR},Q_k^{d,\rR}\ \forall\mbox{ bus }k\in\cA$ (active and reactive, true and reported loads in $\cA$)
\item $P_k^g,Q_k^g\ \forall\mbox{ bus }k\in\cR$ (generation at participating buses)
\item $\forall\mbox{ line } km\in\cE, \ p^\rT_{km}, q^\rT_{km}$, and also $p^\rR_{km}, q^\rR_{km}$ if $km \notin \cL$ (active and reactive, true and reported power flows).
\item $\Delta$  (net change in active power generation)
\end{itemize}

In this formulation, power flows are represented through the  quadratics $S_{km}, S_{mk}$ (see eqs. \eqref{abbreviate}) which appear in the formulation as \eqref{flow1}-\eqref{flow4}.  Note that we include voltage variables but no current variables. However, having solved the above optimization problem, the attacker
reports, for each line $km$ with both ends in $\cA$, a current pair $I^\rR_{km}, I^\rR_{mk}$
computed using the formula
$$\left(^{I^\rR_{km}}_{I^\rR_{mk}}\right) \ = \ Y_{km} \left(^{|V_k^\rR|e^{j \theta^\rR_k}}_{|V_m^\rR|e^{j \theta^\rR_m}}\right),$$
thereby attaining current-voltage consistency.  
Note that if either $k \in \partial\cA$ or $m \in \partial\cA$ the true and reported voltage values are identical -- see Lemma \ref{feas3} below.
\begin{lemma} \label{feas1}
  Consider a feasible solution to problem \eqref{form:attack}.  Let $\rH$ denote
  either $\rT$ or $\rR$ (i.e. true or reported). Then the
  voltages $|V_k^\rH|e^{j \theta_k^\rH}$ for all $k \in \cN$ yield a solution
  to the power flow problem where \begin{itemize}
  \item [(1)] Bus $k$ has load $P_k^{d,\rH} + jQ_k^{d,\rH}$ for $k \in \cA$ and $\hat P_k^{d} + j\hat Q_k^{d}$ if $k \in \cA^C$.
  \item [(2)] Bus $k \in \cG$ has generation $P_k^{g} + jQ_k^{g}$ if $g \in \cR$ and $\hat P_k^{g} + j \hat Q_k^{g}$ if $k \in \cG \setminus \cR$.
  \item [(3)] Line $km$ has power flow $p^\rH_{km} + j q^\rH_{km}$ when $\rH = \rR$ and also when $\rH = \rT$ and $km \notin \cL$.
  \item [(4)] When $\rH = \rR$ (reported data) the solution is fully feasible, i.e. it satisfies voltage, generation, phase angle and power flow limits.
  \item [(5)] When $\rH = \rT$ (true data) the solution satisfies voltage,
    generator and phase angle limits, but only satisfies power flow limits on
    lines $km$ with both $k, m \in \cA^C\cup\partial\cA$. The solution is also consistent with lines in $\cL$ being cut. 
\end{itemize}    
\end{lemma}
\noindent \textit{Proof.} Property (3) follows from constraints \eqref{flow1}-\eqref{flow4}. Hence, (1) and (2) follow from constraints \eqref{rbalwithinA}-\eqref{baloutsideA}. Properties (4)-(5) follow from constraints \eqref{genlimits}-\eqref{linelimits}. \qed

As a corollary to (1)-(2) of Lemma \ref{feas1}, a feasible solution to problem \eqref{form:attack} 
satisfies, exactly, power-injection consistency, i.e. condition (s.4) above.

\begin{lemma} \label{feas2}
  Consider a feasible solution to problem \eqref{form:attack}. The solution
  is consistent with a secondary-response adjustment of active power generator amounting to
  $\Delta$ units.
\end{lemma}
\noindent \textit{Proof.} Follows from constraint \eqref{secondaryresponse}. \qed
\begin{lemma} \label{feas3} Consider a feasible solution to problem \eqref{form:attack}. Then {\bf (a)} the true and reported voltages agree on $\cA^C\cup\partial\cA$.  Further, {\bf (b)} the true and reported currents on a line $km$ are identical if
  $k, m \in \cA^C\cup\partial\cA$.
\end{lemma}  
\noindent \textit{Proof.} {\bf (a)} Follows from constraint \eqref{agreeon},
and {\bf (b)} is a consequence of {\bf (a)}. \qed

\begin{corollary} \label{opt} Suppose we compute a feasible solution
  to problem \eqref{form:attack} whose objective value is strictly greater
  than $(S_{uv}^{max})^2$.  Then the reported solution amounts to an undetectable
  attack that hides an overload on line $uv$.
\end{corollary}
\vspace{-.3in}
\subsection{Computational viability}
Above we have presented a
mathematically correct version of the initial attack problem that would lead to an
(initially) undetectable attack, via problem \eqref{form:attack} which is a nonlinear, nonconvex optimization problem, and thus, in principle, a challenging computational task.  Nevertheless this problem is similar to the standard AC-OPF or PF problem and (at least) a local optimum should be efficiently computable; this expectation is borne out by our experiments.  Strict maximization in
\eqref{form:attack} is \textit{not} required for an attack to be successful (all that is needed is an overload of the line $uv$).

A broader issue concerns the selection of the sets $\cA$ and $\cL$.  This is a
combinatorial problem which is bound to be intractable.  In fact \cite{react}
describes a number of strong NP-hardness results in the DC setting, e.g. given
vectors of phase angles $\theta$ and $\theta'$ it is NP-hard to compute a
set $\cL$ such that $B' \theta' = B \theta$ where $B'$ is the bus susceptance matrix of the network with $\cL$ removed.

Nevertheless, as discussed in the literature, an attacker may be willing to incur significant computational costs in order to compute a successful attack. While it is reasonable to assume that an attacker's 
ability to take physical action or to modify data is limited  (see the discussion in \cite{liureiteretal,poor}), not assuming computational intelligence on the part of
an attacker amounts to a limitation on the part of the defender.

We separate two distinct issues here: first the identification of the set $\cA$, which is done in advance and may be computationally intensive, and second, the solution to \eqref{form:attack} which only requires a few seconds.  Let us assume that the attacker has had (undetected) access to system and sensor data long enough to identify a weak sector of
the transmission system, i.e. the set $\cA$. In this task the attacker
would rely on the fact that typical (time- and day-dependent) load and generation profiles for transmission systems
are statistically predictable with some accuracy. This fact would help the attacker in the computation of a target set $\cA$, perhaps using enumeration, using load estimates in problem \eqref{form:attack}.

Having identified a particular set $\cA$, problem \eqref{form:attack} would be run once again just prior to the attack, now using 
close estimates of the loads obtained from ambient conditions.  Assuming that the attack is perpetrated during a period of slowly changing
loads, and not close in time to a generator redispatch, the attack will likely be sufficiently numerically accurate so as to become difficult to
detect.

\section{Examples}\label{sec:attackex}

\begin{figure}[ht]
  \begin{center}
    \includegraphics[width=0.5\linewidth]{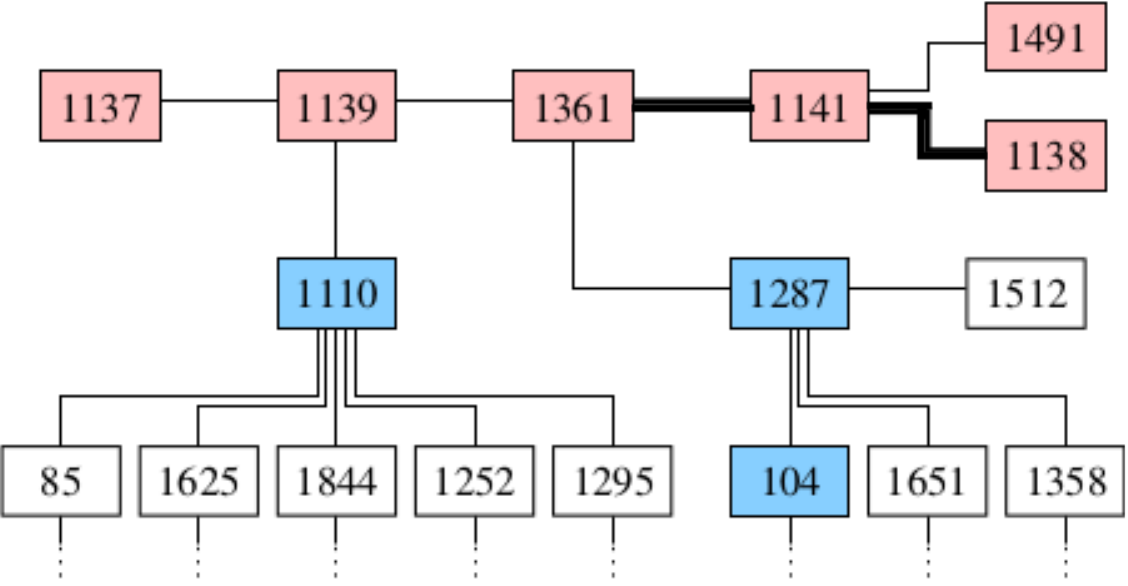}
    \caption{Case2746wp. Attacked zone and its neighborhood. Generators are shown in blue. \label{fig:attack1}}
  \end{center}
    \vspace{-.2in}  
\end{figure}

In the following instance we consider the \texttt{case2746wp} (that has 2746 buses) from the Matpower case library. The adversary attacks the set of buses $\cA= \{$1137, 1138, 1139, 1141, 1361, 1491$\}$ with $\cA - \partial\cA=\{$1137, 1138, 1141, 1491$\}$.  See Figure \ref{fig:attack1}. In this attack the quantity $\Delta$ in \eqref{secondaryresponse} equals $135.09$. We also have $\cL = \emptyset$ (no lines are cut). The set of generators participating in secondary response is $\cR = \{$17, 18, 55, 57, 150, 383, 803, 804, 1996$\}$ with participating factors $\alpha_{k}=1/9$ for all $k\in\cR$.  Note that having a small set
of participating generators makes the attacker's task more challenging, rather than easier, because
it will increase the magnitude of power flows.

\begin{figure}[ht]
  \begin{center}
    \includegraphics[width=0.5\linewidth]{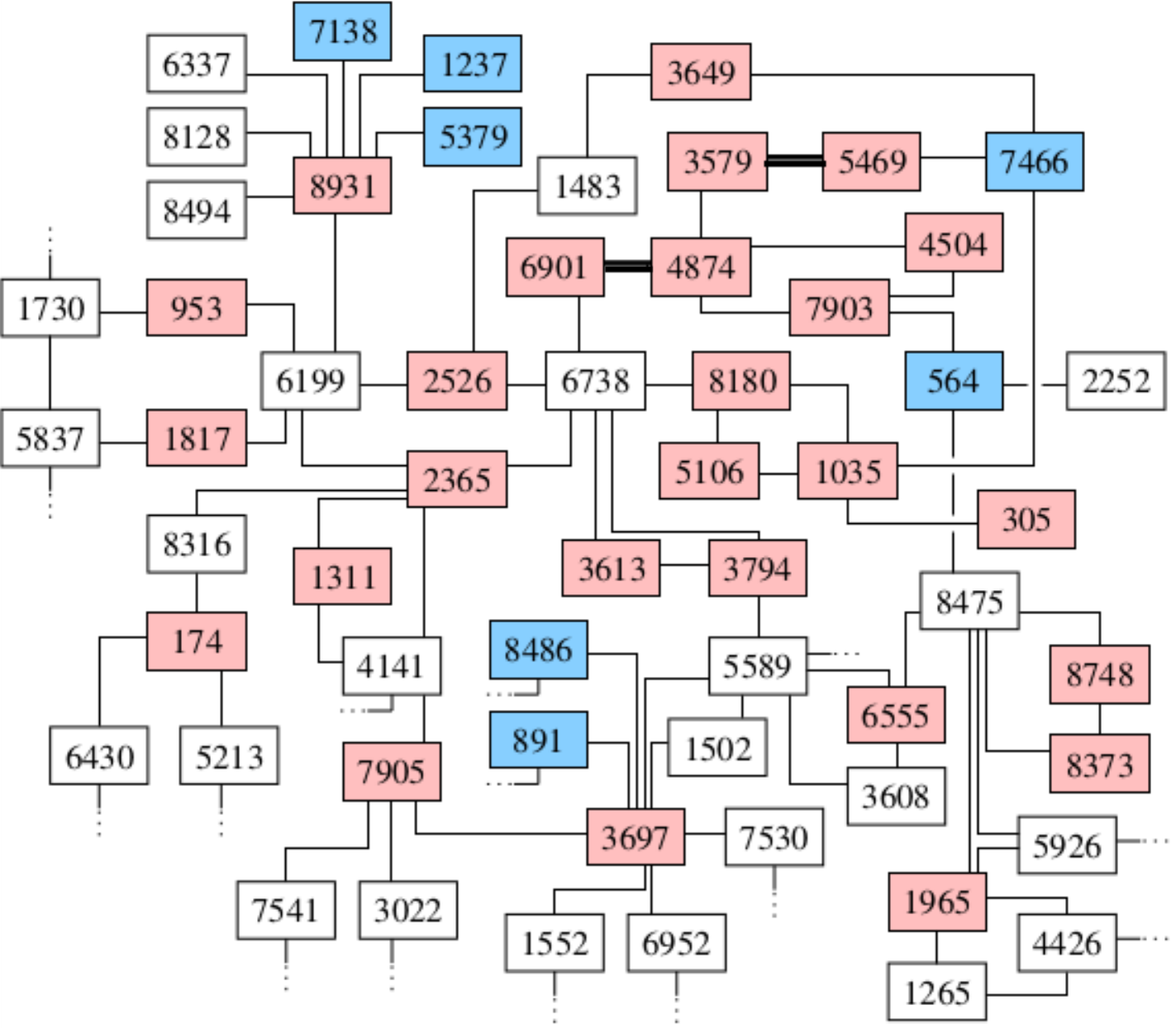}
    \caption{Case1354pegase. Attacked zone and its neighborhood. Generators are shown in blue. \label{fig:attack3}}
  \end{center}
    \vspace{-.2in}  
\end{figure}

Table \ref{tab:attack-flow2} shows the true and reported flow for lines where the solutions differ, with a strong overload on line $(1361, 1141)$ and $(1138, 1141)$. 
As a second example, we consider an ideal attack on the \texttt{case1354pegase} instance of the Matpower library. See Figure~\ref{fig:attack3}. Attack buses are shown in red while generators are shown in blue. The set of participating generators is $\cR = \{$564, 1001, 7466$\}$. 
\begin{table}[ht]
  \begin{center}\vspace{-.1in}
    \caption{True and reported flow at attacked lines. \newline Overloads shown in bold \label{tab:attack-flow2}}
    \scalebox{1}{
      \begin{tabular}{c c | r r | r | c} 
        \multirow{2}{*}{bus $k$} & \multirow{2}{*}{bus $m$} & $p^\rT_{km}$ & $q^\rT_{km}$ & $\|(p_{km}^\rT, q_{km}^\rT) \|$ & \multirow{2}{*}{$S_{km}^{max}$} \\ [0.4ex]
        & & $p^\rR_{km}$ & $q^\rR_{km}$ & $\|(p_{km}^\rR, q_{km}^\rR) \|$ & \\ [0.3ex]
        \hline\hline
        \multirow{2}{*}{1139} & \multirow{2}{*}{1137} & $3.36$ & $2.66$ & $4.29$ \hspace{.3cm} & \multirow{2}{*}{$114.00$} \\
        & & $3.36$ & $2.66$ & $4.28$ \hspace{.3cm} & \\
        \hline
        \multirow{2}{*}{1361} & \multirow{2}{*}{1141} & $229.01$ & $10.49$ & $\bf 229.25$ \hspace{.3cm} & \multirow{2}{*}{$114.00$} \\
        & & $108.51$ & $10.49$ & $109.02$ \hspace{.3cm} & \\
        \hline
        \multirow{2}{*}{1141} & \multirow{2}{*}{1491} & $13.46$ & $2.41$ & $13.68$ \hspace{.3cm} & \multirow{2}{*}{$114.00$} \\
        & & $6.20$ & $2.39$ & $6.64$ \hspace{.3cm} & \\
        \hline
        \multirow{2}{*}{1141} & \multirow{2}{*}{1138} & $209.25$ & $4.44$ & $\bf 209.29$ \hspace{.3cm} & \multirow{2}{*}{$114.00$} \\
        & & $98.06$ & $5.24$ & $98.20$ \hspace{.3cm} & \\
        \hline
      \end{tabular}
    }
  \end{center}
\end{table}
   \vspace{-.2in}  

Table \ref{tab:attack3-flows} shows the true and reported flow for lines where the solutions differ, with line overloads of $66\%$. Detailed solutions of these examples can be found in \cite{solexamples}.

\begin{table}[ht]
\begin{center}
  \vspace{-.1in}
  \caption{True and reported flow for attacked lines. \newline Overloads in shown in bold \label{tab:attack3-flows}}
  \begin{tabular}{c c | r r | r | c } 
    \multirow{2}{*}{bus $k$} & \multirow{2}{*}{bus $m$} & $p^\rT_{km}$ & $q^\rT_{km}$ & $\|(p_{km}^\rT, q_{km}^\rT) \|$ & \multirow{2}{*}{$S_{km}^{max}$}  \\ [0.4ex]
    & & $p^\rR_{km}$ & $q^\rR_{km}$ & $\|(p_{km}^\rR, q_{km}^\rR) \|$ &  \\ [0.3ex] 
    \hline \hline
    \multirow{2}{*}{4874} &  \multirow{2}{*}{4504} & $ -54.84$ & $ -17.66$ & $  57.61$ \hspace{.3cm} &  \multirow{2}{*}{$ 453.00$} \\ 
                                                 & & $ 127.62$ & $  -6.05$ & $ 127.76$ \hspace{.3cm} & \\ 
    \hline
    \multirow{2}{*}{4874} &  \multirow{2}{*}{3579} & $-241.10$ & $  21.74$ & $ 242.08$ \hspace{.3cm} &  \multirow{2}{*}{$ 376.00$} \\ 
                                                 & & $   2.65$ & $  16.17$ & $  16.39$ \hspace{.3cm} & \\ 
    \hline
    \multirow{2}{*}{4874} &  \multirow{2}{*}{3579} & $-254.86$ & $  27.21$ & $ 256.31$ \hspace{.3cm} &  \multirow{2}{*}{$ \infty$} \\ 
                                                 & & $   3.09$ & $  17.07$ & $  17.35$ \hspace{.3cm} & \\ 
    \hline
    \multirow{2}{*}{4874} &  \multirow{2}{*}{7903} & $-202.23$ & $ -56.32$ & $ 209.92$ \hspace{.3cm} &  \multirow{2}{*}{$ 338.00$} \\ 
                                                 & & $ 304.73$ & $ -37.51$ & $ 307.03$ \hspace{.3cm} & \\ 
    \hline
    \multirow{2}{*}{4504} &  \multirow{2}{*}{7903} & $ -56.89$ & $ -13.90$ & $  58.56$ \hspace{.3cm} &  \multirow{2}{*}{$ \infty$} \\ 
                                                 & & $  18.03$ & $ -16.39$ & $  24.37$ \hspace{.3cm} & \\ 
    \hline
    \multirow{2}{*}{3579} &  \multirow{2}{*}{5469} & $-498.23$ & $  34.98$ & $\bf  499.46$ \hspace{.3cm} &  \multirow{2}{*}{$ 491.00$} \\ 
                                                 & & $-120.65$ & $  21.26$ & $ 122.51$ \hspace{.3cm} & \\ 
    \hline
    \multirow{2}{*}{6901} &  \multirow{2}{*}{4874} & $ 953.11$ & $ 236.71$ & $\bf  982.06$ \hspace{.3cm} &  \multirow{2}{*}{$ 591.00$} \\ 
                                                 & & $ 569.95$ & $  88.70$ & $ 576.81$ \hspace{.3cm} & \\ 
    \hline
  \end{tabular}
\end{center}
\vspace{-.1in}
\end{table}

\section{Conclusion}
The possibility of combined physical and data attack on power grids has
gained increased attention.  In principle, an attack that avoids standard
detection methods is possible; we can compute such attacks on large systems
in seconds of CPU time.  This paper focuses on stochastic defense mechanisms
to augment standard detection tools.  Our defenses change the stochastics of
system data in a way that is recognizable by the defender but difficult to
anticipate by the attacker.  In future work we will investigate distributed
versions of our defense mechanisms, as well as protection against more
sophisticated attackers that can also employ distributed resources and
continue the attack over a prolongued period of time.

\section*{Acknowledgment}
This work was supported by DOE award GMLC77; grants
HDTRA1-13-1-0021, ONR N00014-16-1-2889, a DARPA RADICS award, and a DARPA Lagrange award.\\
\vspace{-.25in}
\bibliographystyle{IEEEtran}
\bibliography{sdattack}
\tiny Sat.Aug..3.172842.2019@littleboy
\end{document}